\documentclass[final]{siamltex}
\newtheorem{remark}{Remark}[section]
\newtheorem{example}{Example}[section]
\usepackage{amsmath,tikz,amssymb,graphics,times,multicol}
\usepackage[]{hyperref}

\title{High-order S-Lemma with application to stability of a class of switched nonlinear systems\thanks{
This work is supported by National Natural Science Foundation
of China (No. 61174047), Program for New Century
Excellent Talents in University of Ministry of Education
of China and Basic Research Foundation of Northwestern
Polytechnical University (No. JC201230).
The first author is supported
by the China Scholarship Council.
A short version of this paper
was presented in the 31st Chinese Control Conference, Hefei, China, July 25-27,
2012.}}

\author{Kuize Zhang\thanks{College of Automation,
        Harbin Engineering University, Harbin 150001, PR China
	and Department of Mathematics, University of Turku,
	FIN-20014, Turku, Finland
	({\tt zkz0017@163.com}).}
        \and Lijun Zhang\thanks{School of Marine Technology, Northwestern Polytechnical
	University, Xi'an 710072, PR China and
    College of Automation,
        Harbin Engineering University, Harbin 150001, PR China ({\tt zhanglj7385@nwpu.edu.cn}).}
	\and
	Fuchun Sun\thanks{State Key Laboratory of Intelligent Technology and
	Systems, Tsinghua University, Beijing, 100084, PR China
	({\tt fcsun@mail.tsinghua.edu.cn}).}	}

\begin{document}

\maketitle

\begin{abstract}
This paper extends some results on the S-Lemma proposed by Yakubovich
and uses the improved results to investigate
the asymptotic stability of a class of switched
nonlinear systems.

Firstly, the strict S-Lemma is extended from
quadratic forms to
homogeneous functions with respect to any dilation,
where the improved
S-Lemma is named the strict homogeneous S-Lemma (the
SHS-Lemma for short). In detail, this paper indicates that the strict S-Lemma
does not necessarily hold for homogeneous functions that are not quadratic forms,
and proposes a necessary
and sufficient condition under which the SHS-Lemma holds.

It is well known that a switched linear system with two sub-systems
admits a Lyapunov function with homogeneous derivative (LFHD
for short), if and only if
it has a convex combination of the vector fields of its
two sub-systems that admits a LFHD. In this paper, it is shown that this conclusion
does not necessarily hold for a general switched nonlinear system
with two sub-systems, and gives a necessary
and sufficient condition under which the conclusion holds for a general
switched nonlinear system with two sub-systems.
It is also shown that for a switched nonlinear system
with three or more sub-systems, the ``if'' part holds, but the ``only if'' part
may not.

At last, the S-Lemma
is extended from quadratic polynomials to polynomials of degree more than
$2$ under some mild conditions, and the improved results are called
the homogeneous S-Lemma (the HS-Lemma for short) and the non-homogeneous
S-Lemma (the NHS-Lemma for short), respectively.

Besides, some examples and counterexamples are
given to illustrate the main results.
\end{abstract}



\pagestyle{myheadings}
\thispagestyle{plain}
\markboth{KUIZE ZHANG AND LIJUN ZHANG AND FUCHUN SUN}{HIGH-ORDER S-LEMMA}

\section{Introduction and Preliminaries}\label{intro}

\subsection{S-Lemma}


The S-Lemma, firstly proposed by Yakubovich \cite{Yakubovich:71},
characterizes when a quadratic function is copositive with
another quadratic function. The basic idea of this widely
used method comes from control theory but it has important
consequences in quadratic and semi-definite optimization, convex
geometry, and linear algebra as well  \cite{Polik:07,
Wallin:04}.

A real-valued function $f:R^n\to R$ is said to be copositive with a
real-valued function $g:R^n\to R$ if $g(x) \ge 0$ implies $f(x) \ge 0$.
Furthermore, $f$ is said to be strictly copositive with $g$ if $f$ is
copositive with $g$, and $g(x)\ge0$ and $x\ne0$ imply $f(x)>0$.

\begin{theorem}[S-Lemma, \cite{Yakubovich:71}]\label{thm1}
  Let  $f,g:R^n  \to R$  be quadratic functions such that  $g(\bar x) > 0$
 for some $\bar x \in R^n $. Then  $f$ is copositive with $g$  if and only if there exists   $\xi  \ge 0$
 such that $f(x) - \xi g(x) \ge 0$ for all  $x \in R^n $.
\end{theorem}

\begin{theorem}[strict S-Lemma]\label{thm14}
  Let  $f,g:R^n  \to R$  be quadratic forms. Then  $f$ is strictly
  copositive with $g$  if and only if there exists   $\xi > 0$
  \footnote{In the original version, $\xi\ge0$. However, $f$ is strictly
  copositive with $g$ implies $f$ and $g$ have no common zero point except
  $0\in R^n$. Then by Theorem \ref{thm2}, a positive real number $\xi$ can be
  found.}
 such that $f(x) - \xi g(x) > 0$ for all  nonzero $x \in R^n $.
\end{theorem}

Theorem \ref{thm1} and Theorem \ref{thm14} were firstly obtained
based on the following Theorem \ref{thm2} given in
\cite{Dines:41} via the separation theorem for convex sets.

\begin{theorem}[\cite{Dines:41}]\label{thm2}
  Let  $f,g:R^n  \to R$ be quadratic forms.
  Then the set\\ $\left\{ {\left( {f(x),g(x)} \right):x \in R^n } \right\}$ is convex.
Particularly, if $f$ and $g$ have no common zero point except $0\in R^n$,
then the set $\left\{ {\left( {f(x),g(x)} \right):x \in R^n } \right\}$ is closed
as well as convex, and is either the entire $xy$-plane or an
angular sector of angle less than $\pi$.
\end{theorem}

Yakubovich \cite{Yakubovich:71} gave an example indicating the set
$\{(f(x),g_1 (x),g_2 (x)):x\in R^n\}$ is not convex, which
indicates neither Theorem \ref{thm1} nor Theorem \ref{thm14} holds
for three or more quadratic functions. We shall also give an example to
support it (see Example \ref{thm8}) later. Despite the general
non-convexity of the set  $\{(f(x),g_1 (x),g_2 (x)):x\in R^n\}$,
one can impose additional conditions on quadratic functions
$f(x),g_1(x),\cdots,g_m(x)$ to make the set
$\{(f(x),g_1 (x),\cdots,g_m (x)):x\in R^n\}$ be convex. There are many such
extensions with applications to control theory (linear systems)
\cite{Polik:07,Wallin:04,
Fradkov:73,Yakubovich:92,Uhlig:79,Derinkuyu:06,Li:11}. However, the case
that these functions are (homogeneous) polynomials that
have degree more than $2$
or even general homogeneous functions has not been studied yet,
which can be used to deal with nonlinear systems.
In this paper, we focus on the latter case.

\subsection{Homogeneous Function and Even (Odd) Function}
In this subsection we introduce some preliminaries related to
homogeneous functions and  even (odd) functions.

Any given $n$-tuple $(r_1,\cdots,r_n)$ with each $r_i$
positive is called a dilation;
the set $\{x\in{
R}^n:(|x_1|^{l/r_1}+\cdots+|x_n|^{l/r_n})^{1/l}=1\}$ denotes the
generalized unit sphere, where
$l>0$. Specially, the set $\{x\in{
R}^n:|x_1|^{2}+\cdots+|x_n|^{2}=1\}$ denotes the unit sphere.
Based on the
concept of dilations, the concept of homogeneous
functions is introduced as follows \cite{Rosier:92,Hong:01}:
\begin{definition}\label{homo}
  A function $f:R^n\rightarrow R$
  is said to be homogeneous of degree
  $k\in{R}$ with respect to the dilation $(r_1,\cdots,r_n)$,
  if
  \begin{equation}\label{def_homo_fun}
  f(\epsilon^{r_1}x_1,\cdots,\epsilon^{r_n}x_n)=\epsilon^kf(x_1,\cdots,x_n)
  \end{equation}
  for
  all $\epsilon>0$, and $x_1,\cdots,x_n\in{R}$.
\end{definition}

It can be easily seen that $f$ is homogeneous of degree $k$ with respect
to the dilation $(r_1,\cdots,r_n)$ if and only if $f$ is homogeneous of
degree $k/r$ with respect to the dilation $(r_1,\cdots,r_n)/r$, where
$r=\min\{r_1,\cdots,r_n\}$.
Without loss of generality, we assume that $r_i\ge1$, $i=1,\cdots,n$
hereinafter.
By Definition \ref{homo}, homogeneous polynomials are analytic and
homogeneous functions
of degree a nonnegative integer with respect to the trivial dilation $(1,\cdots,1)$.

A function $f:R^n\to R$ is called even (odd) if $f(-x)=f(x)(-f(x))$ for
all $x\in R^n$. For example, a homogeneous polynomial of even (odd) degree is an
even (odd) function. However, a homogeneous function is not necessarily a polynomial
or not necessarily an even (odd) function.
For example, the odd and homogeneous function $|x|^{\frac{3}{2}}\mbox{sgn}(x)$
is not a polynomial, where $\mbox{sgn}(\cdot)$ denotes the sign function;
the polynomial $x+y^2$ that is homogeneous
of degree $2$ with respect to the dilation $(2,1)$ is neither
an even (odd) function nor a homogeneous polynomial;
the homogeneous function $x^3+|x|^3$ is neither
a polynomial nor an even (odd) function. 

\subsection{Applications of the Strict S-Lemma to Stability of Switched Linear Systems}

Wicks and Peleties \cite{Wicks:94}
showed that if a switched linear system with two sub-systems
has an asymptotically stable convex combination of its sub-systems, there exists a quadratic
Lyapunov function and  a
computable stabilizing switching law. Feron \cite{Feron:96}
proved the converse is also true
by constructing a quadratically
stable convex combination of the two sub-systems based on two total derivatives
(two quadratic forms) of the existing
quadratic Lyapunov function and
using the strict S-Lemma. These results reveal the difference degree
between linear systems and switched linear systems from the perspective of
stability.
Due to their substantial contributions, these results were quoted widely
and embodied in the monograph \cite{Liberzon:03} on switched systems.
However, these
 results have not been extended to nonlinear cases. This is for reason that it
is difficult for switched nonlinear systems to construct stable convex
combinations of the sub-systems, and the strict
S-Lemma can not be used to deal with derivatives of Lyapunov functions of
higher degrees. Then interesting issues arise: May the
strict S-Lemma be extended to nonlinear
functions of higher degrees? May the above necessary and sufficient condition
for switched linear systems be extended to switched nonlinear systems?

Homogeneous nonlinear systems are a class of nonlinear
systems that have properties similar to linear systems, and many
interesting results of linear systems were extended to
homogeneous nonlinear systems (cf. \cite{Rosier:92,Bacciotti:88,Hermes:95,
Hong:01,Aleksandrov:12,Zhang2012}).  Cheng and Martin
\cite{Cheng:01a} proposed the
concept of the Lyapunov function with homogeneous derivative (LFHD for
short)
and applied it to testify the stability of a class of nonlinear
polynomial systems. A nonlinear system admitting a LFHD is not necessarily
homogenous, but still have some properties of homogeneous systems.
For example, if a nonlinear component-wise homogeneous
polynomial system admits a LFHD (cf. \cite{Cheng:01a}), its global stability
is easily guaranteed.
A nonlinear system admitting a LFHD can be regarded as an
approximation of the center manifolds of a large class of nonlinear
systems. Hence to study such systems is theoretically significant and
interesting. Cheng and Martin \cite{Cheng:01a} also gave methods to construct a
LFHD for a component-wise homogeneous polynomial systems.

In this paper, we use the concept of LFHD to characterize a class of
switched nonlinear systems.




\subsection{Model}

In order to describe this problem clearly, the system considered in this
paper is formulated as
\begin{equation}\label{0}
  \dot x = f_{\sigma (t)} (x) ,\quad x = x(t) \in R^n,
\end{equation}
where $\sigma :[0, + \infty ) \to \Lambda  = \{ 1,2, \cdots ,N\} $
is a piece-wise constant, right continuous function, called the
switching signal, $N$ is an integer no less than $2$, and each $f_i$
is a continuous function of the state $x$. A convex combination of the sub-systems
 of system \eqref{0} denotes the system $\dot x=\sum_{i=1}^{N}{\lambda_if_i(x)}$,
 where $0\le\lambda_1,\lambda_2,\cdots,\lambda_N\le1,$ and $\sum_{i=1}^{N}{\lambda_i}=1$.

Throughout this paper, it is assumed that system \eqref{0} admits a
LFHD.
That is to say, there exists a positive definite and continuously
differentiable function
$V:R^n\rightarrow R$, such that each of $\left.{\dot
V(x)}\right|_{S_i}$ is a continuous, even and homogeneous function of the same
degree with respect to the same dilation,
and
\begin{equation}\label{18}\bigcup\limits_{i = 1}^N {\left\{ {x \in R^n
:\left. {\dot V(x)} \right| _{S_i }  < 0} \right\}}  \supset R^n
\setminus\{ 0\} ,\end{equation} where $S_i$ denotes the $i$-th
sub-system, $\left. {\dot V(x)} \right| _{S_i }$ denotes the
derivative of $V(x)$ along the solution trajectory of $S_i$,
$i\in\Lambda$.

It can be proved that if system \eqref{0} admits a LFHD, then for
any given initial state, there exists a switching law driving the initial
state to the equilibrium point as $t\rightarrow \infty$ \cite{Sun2011}.



Based on the concept of LFHD, the necessary and sufficient conditions
given in \cite{Wicks:94,Feron:96} can be restated as:
If for system \eqref{0}, each $f_i$ is linear and $N=2$, then system
\eqref{0} admits a (quadratic) LFHD if and only if there exits a convex
combination of its two sub-systems that admits a (quadratic) LFHD.
In this paper, we will extend these results to nonlinear
system \eqref{0} with $N=2$,
and show that the necessary one does not hold when $N>2$.

The contributions of the paper include:

\begin{itemize}
  \item We extend the strict S-Lemma to the strict homogeneous S-Lemma
        (the SHS-Lemma for short,
	from the case $f,g$  are quadratic forms to
	homogeneous functions with respect to  any dilation).
	In detail, we indicate that the strict S-Lemma does not necessarily hold for
	homogeneous
	  functions that are not quadratic forms, and give a necessary and
	  sufficient condition under which the SHS-Lemma holds.
  \item We use the SHS-Lemma to give a necessary and sufficient condition
	  under which
        system (\ref{0}) when $N=2$  admits a LFHD
        if and only if there exists a convex combination
        of its sub-systems that admits a LFHD, and show
	the ``if'' part still holds when $N>2$.
  \item A counterexample is given to show that even though system (\ref{0})
	when $N>2$ admits a LFHD,
        there may exist no convex combination of its sub-systems
        that admits a LFHD.
  \item The S-Lemma is extended to polynomials of degree more than $2$ under
	  some mild conditions, and the extended results are called the
	  homogeneous S-Lemma (the HS-Lemma) and the non-homogeneous
	  S-Lemma (the NHS-Lemma), respectively.
\end{itemize}

The remaining part of this paper is organized as follows: Section
\ref{sec2} gives the main results and some examples supporting the
main results. The SHS-Lemma is first shown, then
based on it, the asymptotic stability of switched nonlinear systems
with two sub-systems is analyzed;
a counterexample about switched linear systems with
more than two sub-systems is given;
at last, some non-strict S-Lemmas are shown.
Section \ref{sec3} is a brief conclusion.

\section{Main Results}\label{sec2}

Until now, there have been four approaches to proving the S-Lemma
(cf.  \cite{Dines:41,Yakubovich:71}, \cite{Ben-Tal:01,
Sturm:03}, \cite{Yuan:90} and \cite{Hauser:07}, respectively).
It turns out that the two approaches given in \cite{Ben-Tal:01,
Sturm:03} and \cite{Yuan:90} cannot be generalized to prove
the SHS-Lemma, since for homogeneous polynomials of degree more than $2$,
the positive definiteness cannot only be
determined by their coefficient matrices or the eigenvalues of their
coefficient matrices; the approach given in \cite{Hauser:07} cannot
 either, since unlike quadratic polynomials, graphs of
polynomials of degree greater than $2$ are not necessarily spherically
convex (The concept of spherical convexity is referred to
\cite{Hauser:07}). The most fundamental approach, the approach
given in \cite{Dines:41,Yakubovich:71} can be generalized
to deal with the case that the homogeneous functions are odd functions.
However for the case that the homogeneous functions
are even, it does not work either.
In this paper, we propose a new approach that can be used to deal with
both the two cases and to prove the SHS-Lemma.



\subsection{Strict Homogeneous S-Lemma with Application to
Stability of Switched Nonlinear Systems with Two Sub-systems}

We first prove Theorem \ref{thm12} that is an extension of Theorem
\ref{thm2} to some extent, and then prove the SHS-Lemma (Theorem \ref{thm13})
based on Theorem \ref{thm12}.

\begin{theorem}
	Let $f,g:R^n \to R$ be continuous, homogeneous functions of
	degree $0\le k\in R$ with respect to the same dilation
	$(r_1,\cdots,r_n)$, and assume $f$ and $g$ have no common zero point
	except $0\in R^n$ when $k>0$.
	Then the set $\{(f(x),g(x)):x\in R^n\}:=U$ is
	closed. If $k=0$, the set $U$ is a singleton. Next assume $k>0$.
	If $f$ and $g$ are both odd functions, the set $U$ is convex.
	In detail, the set $U$ either equals $R^2$ , or is a straight line
	passing through the origin. If $f$ and $g$ are both even
	functions, the set $U$ is an angular sector.
	\label{thm12}
\end{theorem}

\begin{figure}
\begin{center}
\includegraphics[width=3in]{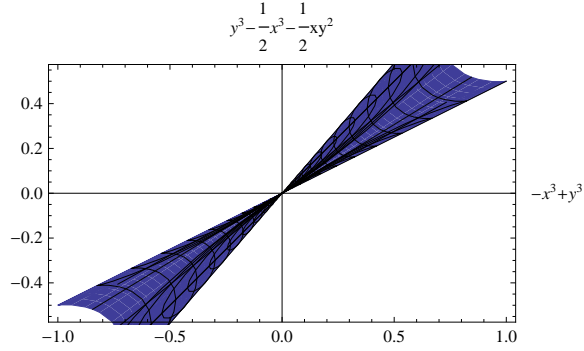}    
\caption{Set $\{(-x^3+y^3,y^3-\frac{1}{2}x^3-\frac{1}{2}xy^2):x^2+y^2\le 1\}$}  
\label{fig10}                                 
\end{center}                              
\end{figure}

\begin{remark}\label{rem1}
	Note that in Theorem \ref{thm12}, the assumption that $f$ and $g$ have
	no common zero point except $0\in R^n$ is crucial. It is because if $f$ and
	$g$ do have a common nonzero zero point, the set $\{(f(x),g(x)):x\in R^n\}$
	may be neither convex nor an angular sector.
	For example, polynomials $f(x,y)=-x^3+y^3$ and $g(x,y)=y^3-\frac{1}{2}x^3-\frac{1}{2}xy^2$ have
	the common nonzero zero point $(1,1)$. And the set
	$\{(-x^3+y^3,y^3-\frac{1}{2}x^3-\frac{1}{2}xy^2):x,y\in R \}:=U$ is neither
    convex nor an angular sector (see Fig. \ref{fig10}).
	This is because $(1,1)$ and $(-1,-\frac{1}{2})$ are both in $U$, but
	$(0,\frac{1}{4})=\frac{1}{2}[(1,1)+(-1,-\frac{1}{2})]$ is not in $U$;
    $(-1,-1)$ and $(1,\frac{1}{2})$ are both in $U$, but
	$(0,-\frac{1}{4})=\frac{1}{2}[(-1,-1)+(1,\frac{1}{2})]$ is not in $U$.

	On the other hand, if $f$ and $g$ have a common nonzero zero point,
	the set $\{(f(x),g(x)):x\in R^n\}$ may be a convex set. For example,
	polynomials $x^2-2xy+y^2$ and $x^2-y^2$ have the common nonzero zero
	point $(1,1)$, but the set $\{(x^2-2xy+y^2,x^2-y^2):x,y\in R\}$ is still
	convex by Theorem \ref{thm2}.
\end{remark}

\begin{proof}[of Theorem \ref{thm12}]
	Let $U$ denote the set $\{(f(x),g(x)):x\in R^n\}$ for short.

	$k=0$:

	Let $\epsilon_m$
	be $1/m$, $m=1,2,\cdots$. We have $\lim_{m\rightarrow\infty}{(\epsilon_m^
	{r_1}x_1,\cdots,\epsilon_m^{r_n}x_n)}=0$ for all $x_1,\cdots,x_n\in
	R$. Further $$f(x_1,\cdots,x_n)=\lim_{m\rightarrow\infty}{\epsilon_m^0
	f(x_1,\cdots,x_n)}=\lim_{m\rightarrow\infty}{f(\epsilon_m
	^{r_1}x_1,\cdots,\epsilon_m^{r_n}x_n)}=f(0)$$
	for all $(x_1,\cdots,x_n)\in R^n$ by the
	continuity and homogeneity of $f$. Similarly $g$ is also constant.
	Hence the set $U$ is a singleton, which is closed.

	$k>0$:

	In this case, $f(0)=g(0)=0$.

	Firstly we prove the set $U$ is closed.

Because $f$ and $g$ are continuous and they have no common zero point
except $0\in R^n$, $(f(x),g(x))/\|(f(x),g(x))\|$ is a continuous function
defined on $R^n\setminus\{0\}$ and maps the unit sphere of $R^n$
onto a compact subset of the unit sphere of $R^2$, where $\|\cdot\|$ is the
Euclidean norm. The compact subset is
also compact in $R^2$, and then closed.
Further by the homogeneity
of $f$ and $g$, the set $U$ is closed.

	Secondly we prove if $u\in U$, then $\lambda u\in U$ for all $\lambda>0$.

	For any given $u\in U$, there exists $z_1=(z_{1_1},\cdots,z_{1_n})\in
	R^n$ such that
	\begin{equation}\label{32}u=(u_f,u_g)=(f(z_1),g(z_1)).
	\end{equation}

	For any given $\lambda>0$, there exists $\bar\epsilon>0$ such that
	$\lambda=\bar\epsilon^k$. Then $$\lambda u=(\bar\epsilon^kf(z_1),
	\bar\epsilon^kg(z_1))=(f(\bar\epsilon^{r_1}z_{1_1},\cdots,
	\bar\epsilon^{r_n}z_{1_n}),g(\bar\epsilon^{r_1}z_{1_1},\cdots,
	\bar\epsilon^{r_n}z_{1_n}))\in U.$$

	When $f$ and $g$ are both odd functions,
    $u\in R$ implies $\lambda u\in U$ for all
	$\lambda\in R$.

	Similar to \eqref{32}, for any given $v\in U$, there exists
	$z_2=(z_{2_1},\cdots,z_{2_n})\in R^n$ such that
	\begin{equation}\label{33}v=(v_f,v_g)=(f(z_2),g(z_2)).
	\end{equation}

	Thirdly we define a closed curve that plays a central role in the
	following proof.
We use $f(\theta)$ and $g(\theta)$ to denote the functions
\begin{equation}
\begin{split}
f( & z_{1_1} |\cos \theta|^
			{r_1}\mbox{sgn}(\cos\theta)  + z_{2_1} |\sin \theta|^{r_1}\mbox{sgn}(\sin\theta),\cdots,\\
		&z_{1_n} |\cos \theta|^{r_n}\mbox{sgn}(\cos\theta)  + z_{2_n} |
		\sin \theta|^{r_n}\mbox{sgn}(\sin\theta))
\end{split}
\end{equation}and
\begin{equation}
  \begin{split}
    g( & z_{1_1} |\cos \theta|^
			{r_1}\mbox{sgn}(\cos\theta)  + z_{2_1} |\sin \theta|^{r_1}\mbox{sgn}(\sin\theta),\cdots,\\
		&z_{1_n} |\cos \theta|^{r_n}\mbox{sgn}(\cos\theta)  + z_{2_n} |
		\sin \theta|^{r_n}\mbox{sgn}(\sin\theta)),
  \end{split}
\end{equation}
respectively for short hereinafter,  where
$\mbox{sgn}(\cdot)$ denotes the sign function.

The function $(f(\theta),g(\theta))$ can be seen as a continuous function
defined over the closed interval $[0,2\pi]$, and $f(\theta)$ and $g(\theta)$
both have period
$2\pi$,
then the curve $\{(f(\theta),g(\theta)):\theta\in[0,2\pi]\}:=\ell$ is a
path-connected, bounded and closed set. And $\{tv_\ell:t\ge0,v_{\ell}\in
\ell\}\subset U$.

Since $f(x)$ and $g(x)$ have no common zero point except $0\in R^n$,
$f(\theta)=g(\theta)=0$ implies $z_{1_i} |\cos \theta|^
{r_i}\mbox{sgn}(\cos\theta)  + z_{2_i} |\sin \theta|^{r_i}\mbox{sgn}
(\sin\theta)=0$, then $z_{1_i}=-z_{2_i}|\tan\theta|^{r_i}\mbox{sgn}(\tan
\theta)$ or $z_{1_i}|\cot\theta|^{r_i}\mbox{sgn}(\cot
\theta)=-z_{2_i}$ for all $i=1,\cdots,n$.
Then $u$ and $v$ are linearly dependent.
Hence the curve $\ell$ does not pass through the origin if $u$ and $v$
are linearly independent.
Similarly, if $f$ and $g$ are both even functions, $u$ and $v$ are linearly
dependent, and either $u_fv_f<0$ or $u_gv_g<0$, the curve $\ell$ is also
path-connected, bounded, closed and does
not pass through the origin either.

At last, we give the conclusion.

Next assume that $f(x)$ and $g(x)$ are both odd functions.

Assume that the set $U$ is not a line passing through
the origin, then there exist linearly independent vectors
$u,v\in U$.
It is easy to get $f(\theta)=-f(\theta+\pi)$ and $g(\theta)=-g(\theta+\pi)$
for all $\theta\in R$. That is, the curve $\ell$ is central symmetric.
Then $\ell$ is homeomorphic to
the unit sphere of $R^2$. Hence $\{tv_{\ell}:
t\ge0,v_{\ell}\in\ell\}= R^2\subset U\subset R^2$. That is, $U=R^2$, and
$U$ is convex.

Next assume that $f$ and $g$ are both even functions.

Assume $U\ne R^2$, that is to say, there exists a vector $u'\in R^2$ such that
$u'\notin U$, then the set $U$ is contained in an angular sector of angle
less than $2\pi$ whose boundary is in $U$ since $U$ is closed.
The boundary of the angular sector is the union of two half lines.
Choose two points $u,v$ in different half lines. Then the corresponding curve
$\ell$ is path-connected, closed and does not pass through the origin.
And furthermore, $\{tv_{\ell}:t\ge0,v_{\ell}\in\ell\}=U$
equals the angular sector.
\end{proof}

\begin{example}\label{exam1}
	We give some examples to illustrate Theorem \ref{thm12}.

	$k$ is odd:

	\begin{enumerate}
		\item $\{(f(x)=x^3,g(x)=x^3):x\in R\}$
			is a straight line passing through the origin..
		\item $\{(f(x_1,x_2)=x_1^3,g(x_1,x_2)=x_2^3):x_1,x_2\in R\}=
			R^2$.
	\end{enumerate}

	$k$ is even:

	In this case, we give some examples to show the angle, denoted by
$\Phi$, of the set $U$ (see the proof of Theorem \ref{thm12}) satisfies $\Phi=\pi$,
$\pi<\Phi<\frac{3}{2}\pi$, $\frac{3}{2}\pi<\Phi<2\pi$ and $\Phi=2\pi$,
respectively. The case $\Phi<\pi$ is seen in Example \ref{thm7}
(see Fig. \ref{fig2}). In each of the following four examples, $f$
and $g$ have no common zero point except $0\in R^2$.

\begin{enumerate}

	\item
$\{(f(x,y)=x^4-y^4-x^2y^2,g(x,y)=-x^4+y^4):x,y\in R\}$ ($\Phi=\pi$):

$(f(1,0),g(1,0))=(1,-1)$,
$(f(0,1),g(0,1))=(-1,1)$ and
$f(x,y)+g(x,y)\le0$ for all $(x,y)\in R^2$ imply
the angle of  $U$ equals $\pi$.

	\item
$\{(f(x,y)=-x^4+y^4-xy^3,g(x,y)=x^4-y^4+x^3y):x,y\in R\}$ ($\pi<\Phi<\frac
{3}{2}\pi$):

$(f(1,-2),g(1,-2))=(23,-17)$,
$(f(2,1),g(2,1))=(-17,23)$ and \\
$(f(3,4),g(3,4))=(-17,-67)$ imply
$(23,-17),(-17,23),(-17,-67)\in U$.
The three points show that the angle of $U$ is greater than $\pi$.
The inequalities $f(x,y)\ge0$ and $g(x,y)\ge0$ have no common solution
shows that the angle of $U$ is less than $\frac{3}{2}\pi$.

	\item
$\{(f(x,y)=x^6-y^6+20x^5y-20x^3y^3,g(x,y)=-x^6+y^6-10xy^5):x,y\in R\}$
($\frac{3}{2}\pi<\Phi<2\pi$):

$(f(0,1),g(0,1))=(-1,1)$,
$(f(2,3),g(2,3))=(-3065,-4195)$,\\
$(f(2,1),g(2,1))=(543,-83)$
and
$(f(-5,-6),g(-5,-6))=(133969,419831)$ imply
$(-1,1),(-3065,-4195),(543,-83),(133969,419831)\in U$. \\
$\langle(543,-83),(133969,419831)\rangle=37899194>0$, where
$\langle\cdot,\cdot\rangle$ denotes the inner product.\\
$f(x,y)=1$ and $g(x,y)=0$ have no common solution.\\
Hence $\frac{3}{2}\pi<\Phi<2\pi$.

	\item
$\{(f(x,y)=x^6-y^6,g(x,y)=-x^6+y^6-x^3y^3):x,y\in R\}$ ($\Phi=2\pi$):

$f(x,y)=a$ and $g(x,y)=b$ have a common solution for all $a,b\in R$.

\end{enumerate}

\end{example}

Based on Theorem \ref{thm12}, we give the following Theorem \ref{thm13}.
We still call it the strict homogeneous S-Lemma.

\begin{theorem}[SHS-Lemma]
	Let $f,g:R^n \to R$ both be continuous, even and homogeneous functions of
	degree $0\le k\in R$ with respect to the same dilation $(r_1,\cdots,r_n)$.
If and only if there exist $a,b\in R$ such that $a^2+b^2>0$ and neither
  $\left\{
  \begin{array}{l}
	  f(x)=a\\
	  g(x)=b
  \end{array}
  \right.
  $ nor $\left\{
  \begin{array}{l}
	  f(x)=-a\\
	  g(x)=-b
  \end{array}
  \right.
  $ have a solution,
the following two items are equivalent:

	$(i)$ $f$ is strictly copositive with $g$;

	$(ii)$ there exists   $\xi > 0$ such that $f(x) - \xi g(x) > 0$  for all  $0 \ne x \in R^n $.
	\label{thm13}
\end{theorem}

\begin{remark}
	Theorem \ref{thm12} shows that if $f$ and $g$ are both homogeneous of
	odd degree and $f$ and $g$ have no nonzero common zero point,
	the set $\{ {(f(x),g(x)):x \in R^n }\}$ is either the whole
	$R^2$ or a straight line passing through the origin. In the former case,
	$(i)$ of Theorem \ref{thm13} cannot hold.
	In the latter case, if $(i)$ of Theorem \ref{thm13} holds,
	there exist $\alpha_1,\alpha_2\in R$ such that
	$\alpha_1\alpha_2<0$ and
	$\alpha_1f(x)+\alpha_2g(x)=0$ for all $x\in R^n$, which indicates
	$(ii)$ of Theorem \ref{thm13}  cannot hold (For example, $f(x)=g(x)
	=x^3:R\to R$.).
	Hence in Theorem \ref{thm13}, we assume that $k$ is even.
\end{remark}

\begin{proof}[of Theorem \ref{thm13}]
	If $k=0$, $f$ and $g$ are both constant functions by Theorem
	\ref{thm12}. Then $(i)$ is obviously equivalent to $(ii)$.

	Next we assume that $k>0$.

  $(ii) \Rightarrow (i)$ holds naturally.

  By Theorem \ref{thm12}, $(i)$ implies
  the set $\{(f(x),g(x)):x\in R^n\}$, denoted by $U$, is an angular sector
  of angle less than $\frac{3}{2}\pi$ and
  $$U\cap\{(r_1,r_2):r_1\le0,r_2\ge0\}=\emptyset.$$

  Next we assume that there exist $a,b\in R$ such that $a^2+b^2>0$ and neither
  $\left\{
  \begin{array}{l}
	  f(x)=a\\
	  g(x)=b
  \end{array}
  \right.
  $ nor $\left\{
  \begin{array}{l}
	  f(x)=-a\\
	  g(x)=-b
  \end{array}
  \right.
  $ have a solution and
  prove $(i) \Rightarrow  (ii)$.

  The foregoing assumption and $(i)$
  imply the angle of $U$ is less than $\pi$. Then
  there exist $\xi_1<0$ and $\xi_2>0$ such that
  $$\xi_1f(x)+\xi_2g(x)<0$$
  for all $0\ne x\in R^n$. Set $\xi=-\xi_2/\xi_1>0$, then $f(x)-\xi g(x)>0$
  for all $x\in R^n$.

  In particular,  when $k=2$, $(i)$ implies the above assumption
  (see Theorem \ref{thm2}).

  Next we assume for all
  $a,b\in R$ such that $a^2+b^2>0$, either
  $\left\{
  \begin{array}{l}
	  f(x)=a\\
	  g(x)=b
  \end{array}
  \right.
  $ or $\left\{
  \begin{array}{l}
	  f(x)=-a\\
	  g(x)=-b
  \end{array}
  \right.
  $ have a solution, which together with $(i)$
  implies the angle of $U$ is no less than $\pi$. Hence $(ii)$ does not hold.

\end{proof}

%

Based on Theorem \ref{thm13}, we give the following Theorem \ref{thm15}.

\begin{theorem}
System (\ref{0}) admits a LFHD $V:R^n\to R$, and
  there exist $a,b\in R$ such that $a^2+b^2>0$ and neither
  $\left\{
  \begin{array}{l}
	  \left.\dot V(x)\right|_{S_1}=a\\
	  \left.\dot V(x)\right|_{S_2}=b
  \end{array}
  \right.
  $ nor $\left\{
  \begin{array}{l}
	  \left.\dot V(x)\right|_{S_1}=-a\\
	  \left.\dot V(x)\right|_{S_2}=-b
  \end{array}
  \right.
  $ have a solution,
  if and only if there exists a convex
combination of its two sub-systems that admits a LFHD when $N=2$ (The ``if
part'' still holds when $N>2$).
	\label{thm15}
\end{theorem}

\begin{proof}
  ``if'':  This part is trivial just like the triviality of the ``if'' part
  of the SHS-Lemma.
  ``only if'': This part is proved by the SHS-Lemma.

  Since $V$ is a LFHD of system \eqref{0} when $N=2$, that is to say,

    $$\bigcup\nolimits_{i = 1}^2 {\left\{ {x \in R^n :\left. {\dot V(x)} \right|_{S_i}
  < 0} \right\}}  \supset R^n \setminus  \left\{ 0
  \right\},$$ then $-\left. {\dot V(x)} \right|_{S_1}$ is strictly copositive
  with $\left. {\dot V(x)} \right|_{S_2}$.

  By Theorem \ref{thm13} and the assumption related to $V$ in Theorem \ref{thm15},
  there exists $\xi  >
  0$ such that $$\left. {\dot V(x)} \right|_{S_1 }  + \xi \left. {\dot V(x)} \right|_{S_2}   <
  0\mbox{ for all }0 \ne x \in R^n.$$
  Take $\lambda _1  = \frac{1}{{1 + \xi }}$, $\lambda _2  = \frac{\xi }{{1 + \xi
  }}$, then $V$ is a LFHD of system
  $\dot x = \lambda _1 f_1(x)  + \lambda _2 f_2(x). $
  \end{proof}

  \begin{remark}
    Theorem \ref{thm15} indicates the existence of an asymptotically stable convex
    combination of the two sub-systems, but it does not show how to
    find the convex combination. Luckily, there are only two sub-systems,
    so we can use Young's inequality to construct the convex
    combination. Example \ref{thm7} illustrates the procedure and
    the case that $\Phi<\pi$ in Theorem \ref{thm12} by showing a switched
    polynomial system and Example \ref{thm16} illustrates the procedure
    by showing a switched non-polynomial system.

    In fact, Theorem \ref{thm15} supplies a method to find a LFHD for
    a switched polynomial system with two sub-systems: $(i)$ Construct its convex
    combination of its sub-systems with coefficients variable parameters;
    $(ii)$ construct a LFHD by using the methods proposed in
    \cite{Cheng:01a}.
  \end{remark}

  \begin{example}\label{thm7}
    Consider the switched polynomial system $S$ with two sub-systems as follows:
    \begin{equation*}
	    \begin{split}
		    &S_1 :\left\{ \begin{array}{l}
 \dot x_1  = 7x_1^3  - 3x_2^3 + 2x_1x_2^2,  \\
 \dot x_2  = 5x_1^3  - 5x_2^3,  \\
 \end{array} \right.\\
&S_2 :\left\{ \begin{array}{l}
 \dot x_1  =  - 5x_1^3 - x_1x_2^2,  \\
 \dot x_2  =  - x_1^3  + x_2^3.  \\
 \end{array} \right.
 \end{split}\end{equation*}

  It is obvious that the origin is the unique equilibrium point for both
  sub-system  $S_1$ and sub-system $S_2$.

  Firstly, we prove the origin is unstable both for sub-system  $S_1$ and for sub-system
  $S_2$.

  For sub-system  $S_1$, choose $V_1 (x) = \frac{1}{4}\left( {5x_1^4  - x_2^4 }
  \right)$. On the line $x_2=0$, $V_1 (x) > 0$  at points arbitrarily close
  to the origin, and $\dot V_1 (x) = 15x_1^6  + 5\left( {2x_1^3  - x_2^3 }
  \right)^2 + 10x_1^4x_2^2
  $ is positive definite. Then by Chetaev's theorem (Theorem 4.3 of \cite{Khalil:07}), the origin is
  unstable.

  For sub-system  $S_2$, choosing $V_2 (x) = \frac{1}{4}\left( {-x_1^4  + x_2^4 }
  \right)$, similarly we have the origin is unstable.

  Secondly, we prove switched system $S$ admits a LFHD.

  Choosing $V(x) = \frac{1}{4}\left( {x_1^4  + x_2^4 } \right)$  that is positive definite,
  then $$ \begin{array}{l}
 \left. {\dot V(x)} \right|_{S_1 }  = 7x_1^6  + 2x_1^3 x_2^3  - 5x_2^6 + 2x_1^4x_2^2, \\
 \left. {\dot V(x)} \right|_{S_2 }  =  - 5x_1^6  - x_1^3 x_2^3  + x_2^6 - x_1^4x_2^2. \\
 \end{array} $$

 By Young's inequality, we have (see Fig. \ref{fig4})
  \begin{equation}\label{9}
	  \begin{split}
      &\left\{(x_1 ,x_2 ): \left. \dot V(x) \right|_{S_1 } < 0\right\}
      \cup  \left\{(x_1 ,x_2 ): \left. \dot V(x) \right|_{S_2 } < 0\right\}\\
      \supset& R ^2 \setminus \{ (0,0)\}.
      \end{split}
  \end{equation}

  The procedure is as follows:  By Young's inequality, we have
 \begin{equation*}
	 \begin{split}
 \left. {\dot V(x)} \right|_{S_1 } & \le \frac{{25}}{3}x_1^6  + 2x_1^3 x_2^3  - \frac{{13}}{3}x_2^6  \\
 &=  - \frac{{13}}{3}\left( {x_2^3  - \frac{{3 - \sqrt {334} }}{{13}}x_1^3 } \right)\left( {x_2^3  - \frac{{3 + \sqrt {334} }}{{13}}x_1^3 } \right), \\
 \left. {\dot V(x)} \right|_{S_2 } & \le  - 5x_1^6  - x_1^3 x_2^3  + x_2^6  \\
= &\left( {x_2^3  - \frac{{1 - \sqrt {21} }}{2}x_1^3 } \right)\left( {x_2^3  - \frac{{1 + \sqrt {21} }}{2}x_1^3 } \right),
 \end{split}
 \end{equation*}
  and then
 \begin{equation*}
	 \begin{split}
  &\left\{ {(x_1 ,x_2 ): \left.{\dot V(x)} \right|_{S_1 } < 0} \right\} \\
  \supset &\left\{ {(x_1 ,x_2 ):\frac{{25}}{3}x_1^6  + 2x_1^3 x_2^3  - \frac{{13}}{3}x_2^6  < 0}
  \right\}\\
  =&\left\{ (x_1 ,x_2 ): {x_2  - \left(\frac{{3 - \sqrt {334} }}{{13}}\right)^{1/3}x_1 } >
  0,
   {x_2  - \left(\frac{{3 + \sqrt {334} }}{{13}}\right)^{1/3}x_1 } > 0 \right\} \\
  &\cup \left\{ (x_1 ,x_2 ): {x_2  - \left(\frac{{3 - \sqrt {334} }}{{13}}\right)^{1/3}x_1 }
  <0,
  {x_2  - \left(\frac{{3 + \sqrt {334} }}{{13}}\right)^{1/3}x_1 } < 0 \right\} ,
  \end{split}
 \end{equation*}
 \begin{equation*}
	 \begin{split}
 &\left\{ {(x_1 ,x_2 ): \left.{\dot V(x)} \right|_{S_2 } < 0} \right\} \\
  \supset &\left\{ {(x_1 ,x_2 ): - 5x_1^6  - x_1^3 x_2^3  + x_2^6  < 0}
  \right\}\\
  = &\left\{ (x_1 ,x_2 ): {x_2  - \left(\frac{{1 - \sqrt {21} }}{{2}}\right)^{1/3}x_1 } >
  0,
  {x_2  - \left(\frac{{1 + \sqrt {21} }}{{2}}\right)^{1/3}x_1 } < 0 \right\} \\
 & \cup \left\{ (x_1 ,x_2 ): {x_2  - \left(\frac{{1 - \sqrt {21} }}{{2}}\right)^{1/3}x_1 }
  <0,
   {x_2  - \left(\frac{{1 + \sqrt {21} }}{{2}}\right)^{1/3}x_1 } > 0 \right\} .
  \end{split}
 \end{equation*}

  As $ {\frac{{1 + \sqrt {21} }}{2}}   > {\frac{{3 + \sqrt {334} }}{13}}  >
  0 > {\frac{{3 - \sqrt {334} }}{13}} > {\frac{{1 - \sqrt {21}
  }}{2}},$ we have
    \begin{equation*}
	    \begin{split}
      &\left\{(x_1 ,x_2 ): \left. \dot V(x) \right|_{S_1 } < 0\right\}
      \cup  \left\{(x_1 ,x_2 ): \left. \dot V(x) \right|_{S_2 } <   0\right\}\\
      \supset&\left\{ {(x_1 ,x_2 ):{\frac{{25}}{3}x_1^6  + 2x_1^3 x_2^3  - \frac{{13}}{3}x_2^6 } < 0}
      \right\}
      \cup \left\{ {(x_1 ,x_2 ): { - 5x_1^6  - x_1^3 x_2^3  + x_2^6 } < 0}\right\} \\
      \supset& R ^2 \setminus \{ (0,0)\}.
      \end{split}
    \end{equation*}

    Thirdly, we prove the LFHD $V$ satisfies the assumption in Theorem \ref{thm15}.

    $\left.\dot V(x)\right|_{S_1}=2$ and $\left.\dot V(x)\right|_{S_2}=-1$
    imply $x_1^6+x_2^6=0$, then $x_1=x_2=0$. That is to say, they have no
    common solution.

    $\left.\dot V(x)\right|_{S_1}=-2$ and $\left.\dot V(x)\right|_{S_2}=1$
    also imply $x_1^6+x_2^6=0$, then $x_1=x_2=0$, which also means they
    have no common solution.

\begin{figure}
\begin{center}
\includegraphics[height=10cm]{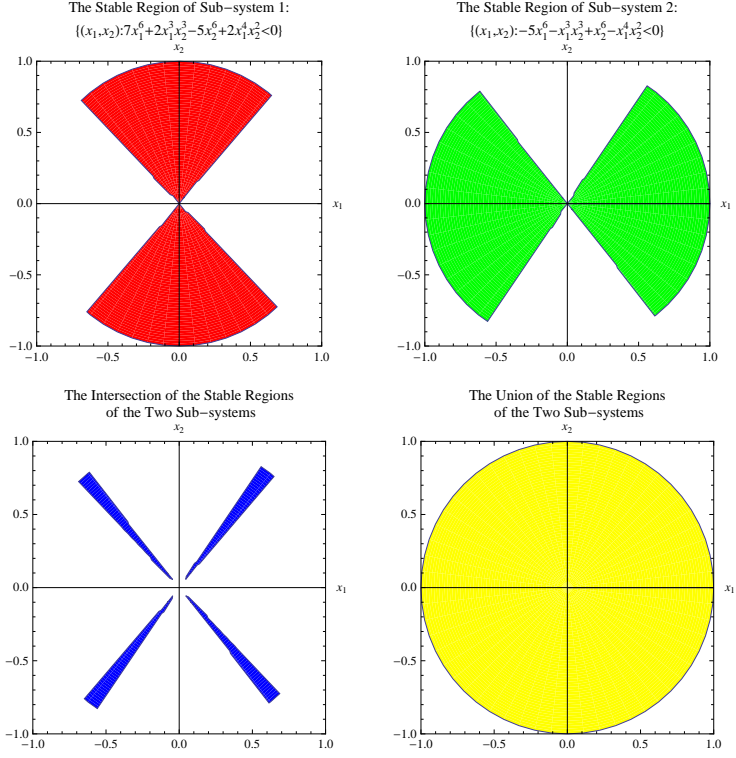}    
\caption{The stable regions of Example \ref{thm7} on the unit disk}  
\label{fig4}                                 
\end{center}                                 
\end{figure}

  To illustrate Theorem \ref{thm12}, the best we can do is to picture the set $$
 \left\{ {\left( {\left. \dot V(x) \right|_{S_1 }, \left. \dot V(x) \right|_{S_2 } } \right):x_1 ,x_2  \in R} \right\}:=U(see Fig. \ref{fig2}).
 $$ From Fig. \ref{fig2} we see that $U$ is an angular sector of angle less than  $\pi$.

  At last, we construct a convex combination of sub-system  $S_1$ and sub-system
  $S_2$ that admits a LFHD.

  Let  $0<\lambda<1$, then a convex combination of sub-system  $S_1$ and sub-system   $S_2$, $\lambda S_1  + (1 - \lambda )S_2
  $,
  is formulated as follows: \begin{equation}\label{12}\left\{ \begin{array}{l}
 \dot x_1  = (12\lambda  - 5)x_1^3  - 3\lambda x_2^3 + (3\lambda - 1) x_1x_2^2, \\
 \dot x_2  = (6\lambda  - 1)x_1^3  + (1 - 6\lambda )x_2^3 . \\
 \end{array} \right.\end{equation}

  We might as well take  $V(x) = \frac{1}{4}\left( {x_1^4  + x_2^4 } \right)$, a positive definite function, then we have
  \begin{equation}\label{10}
	  \begin{split}
    \dot V(x) =& (12\lambda  - 5)x_1^6  + (3\lambda  - 1)x_1^3 x_2^3\\
    &+ (1 - 6\lambda)x_2^6 + (3\lambda - 1) x_1^4x_2^2.
    \end{split}
  \end{equation}

  Now we try to find a $\lambda\in(0,1)$ such that (\ref{10}) is
  negative definite.

  If $3\lambda - 1 \ge 0$, by Young's inequality, we get
  $$\dot V(x) \le \left(\frac{31}{2} \lambda - \frac{37}{6}\right) x_1^6 +
  \left(\frac{1}{6} - \frac{7}{2} \lambda \right) x_2^6.$$

  Let $\frac{31}{2} \lambda - \frac{37}{6}<0$ and $\frac{1}{6} - \frac{7}{2}
  \lambda<0$, together with $3\lambda - 1 \ge 0$, we get $\frac{1}{3} \le \lambda < \frac{37}{93}$.

  If $3\lambda - 1 \le 0$, by Young's inequality, we get
  $$\dot V(x) \le \left(9 \lambda - \frac{9}{2}\right) x_1^6 +
  \left(\frac{3}{2} - \frac{15}{2} \lambda \right) x_2^6+
  (3\lambda-1)x_1^4x_2^2.$$

  Let $9 \lambda - \frac{9}{2}<0$, $\frac{3}{2} - \frac{15}{2}
  \lambda<0$ and $3\lambda-1<0$,
   we get $\frac{1}{5} < \lambda < \frac{1}{3}$.

  Hence, if $\frac{1}{5} < \lambda < \frac{37}{93}$, system (\ref{12})
  admits a LFHD.
  \end{example}

\begin{figure}
\begin{center}
\includegraphics[]{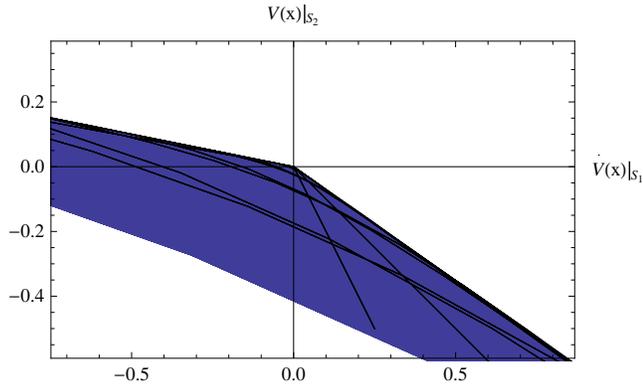}    
\caption{The set $
 \left\{ {\left( {\left. \dot V(x) \right|_{S_1 }, \left. \dot V(x) \right|_{S_2 } }
 \right):x_1^2+x_2^2 \le 1} \right\}$ of Example \ref{thm7}}  
\label{fig2}                                 
\end{center}                                 
\end{figure}

\begin{example}\label{thm16}
    Consider the the following switched system $S$ with two sub-systems as follows:
    \begin{equation*}
	    \begin{split}
		    S_1 :\left\{ \begin{array}{l}
 \dot x_1  = -4x_1 ,\\
 \dot x_2  = 4x_1^{\frac{2}{3}}x_2  + 4x_2^3,  \\
 \end{array} \right.\quad
 S_2 :\left\{ \begin{array}{l}
	\dot x_1  =  2x_1 + x_1^{\frac{1}{3}}x_2^2,  \\
 \dot x_2  =  -8 x_2^3.  \\
 \end{array} \right.
 \end{split}\end{equation*}

  It is obvious that the origin is the unique equilibrium point for both
  sub-system  $S_1$ and sub-system $S_2$.

  Firstly, we prove the origin is unstable both for sub-system  $S_1$ and for sub-system
  $S_2$.

  For sub-system  $S_1$, choose $V_1 (x) = -3 x_1^{\frac{4}{3}}  + x_2^2
  $. On the line $x_1=0$, $V_1 (x) > 0$  at points arbitrarily close
  to the origin, and $\dot V_1 (x) = (4x_1^{\frac{2}{3}} + x_2^2)^2 + 7x_2^4
  $ is positive definite. Then by Chetaev's theorem (Theorem 4.3 of \cite{Khalil:07}), the origin is
  unstable.

  For sub-system  $S_2$, choosing $V_2 (x) = 3 x_1^{\frac{4}{3}}  - x_2^2
  $, then  $\dot V_2 (x) = 4x_1^{\frac{4}{3}} +
  (2x_1^{\frac{2}{3}} + x_2^2)^2 + 15x_2^4$,
  similarly we have the origin is unstable.

  Secondly, we prove switched system $S$ admits a LFHD.

  Choosing $V(x) = 3 x_1^{\frac{4}{3}}  + x_2^2$
  that is positive definite,
  then $$ \begin{array}{l}
	  \left. {\dot V(x)} \right|_{S_1 }  = -16x_1^{\frac{4}{3}}  + 8x_1^
	  {\frac{2}{3}} x_2^2  + 8x_2^4 ,\\
 \left. {\dot V(x)} \right|_{S_2 }  =   8x_1^{\frac{4}{3}}  + 4x_1^
	  {\frac{2}{3}} x_2^2  - 16x_2^4 ,\\
 \end{array} $$which are both homogeneous functions of degree $4$ with respect
 to the  dilation
 $(3,1)$.

 By Young's inequality, we have
  \begin{equation}
	  \begin{split}
      &\left\{(x_1 ,x_2 ): \left. \dot V(x) \right|_{S_1 } < 0\right\}
      \cup  \left\{(x_1 ,x_2 ): \left. \dot V(x) \right|_{S_2 } < 0\right\}\\
      \supset& R ^2 \setminus \{ (0,0)\}(\mbox{see Fig. }\ref{fig6}).
      \end{split}
  \end{equation}

  Thirdly, we prove the LFHD $V$ satisfies the assumption in Theorem \ref{thm15}.

    $\left.\dot V(x)\right|_{S_1}=1$ and $\left.\dot V(x)\right|_{S_2}=-1$
    imply $2x_1^{\frac{4}{3}}-3x_1^{\frac{2}{3}}x_2^2+2x_2^4=0$, which has no
    solution. That is to say, they have no common solution.

    $\left.\dot V(x)\right|_{S_1}=-1$ and $\left.\dot V(x)\right|_{S_2}=1$ also
    imply $2x_1^{\frac{4}{3}}-3x_1^{\frac{2}{3}}x_2^2+2x_2^4=0$, then
    they have no common solution.

\begin{figure}
\begin{center}
\includegraphics[height=10cm]{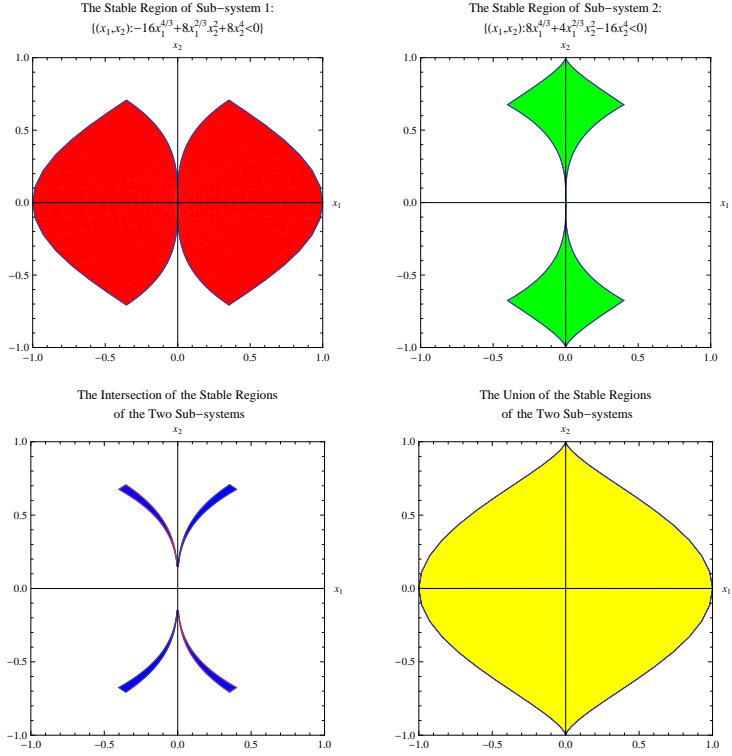}    
\caption{The stable regions of Example \ref{thm16} on the generalized unit disk
$\{(x_1,x_2):x_1^{\frac{2}{3}}+x_2^2\le1\}$}  
\label{fig6}                                 
\end{center}                                 
\end{figure}

  At last, we construct a convex combination of sub-system  $S_1$ and sub-system
  $S_2$ that admits a LFHD.

  Let  $0<\lambda<1$, then a convex combination of sub-system  $S_1$ and sub-system   $S_2$, $\lambda S_1  + (1 - \lambda )S_2
  $,
  is formulated as follows:
  \begin{equation}\label{41}
	  \left\{ \begin{array}{l}
		  \dot x_1  = (2-6\lambda )x_1 + (1-\lambda) x_1^{\frac{1}{3}}x_2^2, \\
		  \dot x_2  = 4\lambda  x_1^{\frac{2}{3}}x_2  + (12\lambda-8 )x_2^3 . \\
 \end{array} \right.\
 \end{equation}

  We might as well take  $V(x) = 3 x_1^{\frac{4}{3}}  + x_2^2$, a positive definite function, then we have
  \begin{equation}\label{40}
	  \begin{split}
    \dot V(x) = &8(1-3\lambda)x_1^{\frac{4}{3}}  + 4(1+\lambda)x_1^
	  {\frac{2}{3}} x_2^2  + 8(3\lambda-2)x_2^4 ,\\
	      \le & (10-22\lambda)x_1^{\frac{4}{3}}  +
	  (26\lambda-14)x_2^4 \\
  	\end{split}
  \end{equation}by Young's inequality.

  Now we try to find a $\lambda\in(0,1)$ such that (\ref{40}) is
  negative definite.

  Let $10-22\lambda<0$ and $26\lambda-14<0$,
   we get $\frac{5}{11} < \lambda < \frac{7}{13}$.

  Hence, if $\frac{5}{11} < \lambda < \frac{7}{13}$, system (\ref{41})
  admits a LFHD.
  \end{example}

  Next we give a direct corollary of Theorem \ref{thm12}.
  We use a generalization of the basic idea in \cite{Dines:41} to
  give an interesting proof that is only suitable for homogeneous
  polynomials of odd degree.

\begin{corollary}\label{thm3}
  Let $f,g:R^n  \to R$ be homogeneous polynomials of degree $k\ge1$,
  and assume that $f$ and
  $g$ have no common zero point except $0\in R^n$. Then the set
  $\{(f(x),g(x)):x\in R^n\}$, denoted by $U$, is closed.
  If $k$ is odd, the set $U$ is convex. In detail, the set $U$ either equals $R^2$,
  or is a straight line passing through
  the origin. If $k$ is even, the set $U$ is
  an angular sector.
\end{corollary}

\begin{proof}

  Let $U$ denote the set $\left\{ {\left( {f(x),g(x)} \right):x \in R^n } \right\}$ for
  short.

  If $A\in U$, then each point in
  the ray starting at the origin and passing through $A$ is in the set $U$ by the
  homogeneity of $f$ and $g$. Hereinafter, we assume that $f$ and $g$ have no
  common zero point except $0\in R^n$ and $k$ is odd.

  Next we prove the set $U$ is a convex set.

  To this end, we only need to prove that for any $u,v \in U$,
  $\lambda u + (1 - \lambda )v \in U$ for all $\lambda  \in [0,1]$.

  There exist $z_1 ,z_2  \in R^n $ such that $$u_f  = f(z_1 ),\mbox{ } u_g  = g(z_1 ),\mbox{ }v_f  = f(z_2 )\mbox{ and } v_g  = g(z_2
  ),$$ where $u = (u_f ,u_g )$ and $v = (v_f ,v_g )$.

  By the homogeneity of $f$ and $g$, if $u$ and $v$ are linearly dependent, then $\lambda u + (1 - \lambda )v \in U$
  for all $\lambda  \in [0,1]$.

  Without loss of generality, we assume that $u$ and $v$ are linearly independent and
  \begin{equation}
  \label{d}
  u_g v_f  - u_f v_g : = d >  0.
  \end{equation}

  Below we try to find a vector $z \in R^n $ such that
  \begin{equation}\label{3}
    \left( {f(z),g(z)} \right) = \lambda u + (1 - \lambda )v
  \end{equation}
  for some $\lambda  \in (0,1)$.

  We make the following ansatz
  $z = \rho \left( {z_1 \cos \theta  + z_2 \sin \theta }
  \right)$, where $\rho$ and $\theta$ are real variables.

  Substitute $z = \rho \left( {z_1 \cos \theta  + z_2 \sin \theta }
  \right)$ into (\ref{3}), we get
  \begin{equation}\label{4}
    \left\{ \begin{array}{l}
 \rho ^{k} f(z_1 \cos \theta  + z_2 \sin \theta ) = \lambda u_f  + (1 - \lambda )v_f , \\
 \rho ^{k} g(z_1 \cos \theta  + z_2 \sin \theta ) = \lambda u_g  + (1 - \lambda )v_g . \\
 \end{array} \right.
  \end{equation}

  Hereinafter, we use $f(\theta)$ and $g(\theta)$ to denote
  $f(z_1\cos\theta+z_2\sin\theta)$ and $g(z_1\cos\theta+z_2\sin\theta)$,
  respectively for short. Then there exists no $\theta'$ such that
  $f(\theta')=g(\theta')=0$. This is because if there does exist
  $\theta'$ such that $f(\theta')=g(\theta')=0$, then $z_1\cos\theta'+z_2\sin\theta'=0$, that is to say, $z_1$ and $z_2$ are linearly dependent; furthermore,
  $u$ and $v$ are linearly
  dependent, which is a contradiction.
  (\ref{4}) shows that
  $\rho^k=d/T(\theta)$ and $\lambda  = S(\theta )/T(\theta )$, where
  \begin{equation}\label{5}
 \begin{split}
 T(\theta ) = &f(\theta )(u_g  - v_g )
 - g( \theta )(u_f  - v_f ), \\
 S(\theta ) = &g( \theta )v_f
 - f( \theta )v_g . \\
 \end{split}
  \end{equation}

  Denote $ S(\theta )/T(\theta ):=\Lambda (\theta )$, a function of
  $\theta$ having period $\pi$. Then
  we need to prove
  \begin{equation}\label{judgement}
	  \Lambda
  \left( {[0,2\pi ] \cap \left\{ {\theta :T(\theta ) > 0} \right\}} \right)
  \supset[0,1].
  \end{equation}

  It is easy to get $S(0) = T(0) = T\left( {\frac{\pi }{2}} \right) =
  d$ and $S\left( {\frac{\pi }{2}} \right) = 0$. So, $\Lambda (0) =
  1$ and $\Lambda \left( {\frac{\pi }{2}} \right) = 0$. Since $f$
  and $g$ are homogeneous polynomials of $\cos \theta $ and
  $\sin \theta$, $T$ and $S$ can be expressed as $$\left\{ \begin{array}{l}
 T(\theta ) = \sum\limits_{i = 0}^{k} {\alpha _i \cos ^i \theta \sin ^{k - i} \theta } , \\
 S(\theta ) = \sum\limits_{i = 0}^{k} {\beta _i \cos ^i \theta \sin ^{k - i} \theta } . \\
 \end{array} \right.$$

  It is obvious that $S(0) = T(0) = T\left( {\frac{\pi }{2}} \right) = d$ and $S\left( {\frac{\pi }{2}} \right) =
  0$, then $\alpha _0  = \alpha _{k}  = \beta _{k}  = d$ and $\beta _0  =
  0$. Hence, $$\left\{ \begin{array}{l}
 T(\theta ) = d\left( {\cos ^{k} \theta  + \sin ^{k} \theta } \right) + \sum\limits_{i = 1}^{k - 1} {\alpha _i \cos ^i \theta \sin ^{k - i} \theta } , \\
 S(\theta ) = d\cos ^{k} \theta  + \sum\limits_{i = 1}^{k - 1} {\beta _i \cos ^i \theta \sin ^{k - i} \theta } . \\
 \end{array} \right.$$

  Notice that $T$ and $S$ are both continuous functions of $\theta$, if $T(\theta ) \ne
  0$ for all $\theta  \in [0,\frac{\pi}{2} ]$ or $[\frac{\pi}{2},\pi] $,
  then $\Lambda (\theta )$ is also a continuous function defined on the interval  $[0,\frac{\pi}{2} ]$ or $[\frac{\pi}{2},\pi] $.
   So
  \begin{equation}\label{6}
    \Lambda \left( {[0,\pi ]} \right) \supset [0,1].
  \end{equation}

  Next we assume that $T$ has a zero point.

  We claim that $T$ and $S$
   have no common zero point.
  If there exists $\hat \theta  \in \left( {0,\frac{\pi }{2}} \right) \cup \left( {\frac{\pi }{2},\pi }
  \right)$ such that $T(\hat \theta ) = S(\hat \theta ) = 0$, then
  \begin{equation}\label{7}
    \left[ {\begin{array}{*{20}c}
   { - v_g } & {v_f }  \\
   {u_g } & { - u_f }  \\
\end{array}} \right]\left[ {\begin{array}{*{20}c}
   {f(z_1 \cos \hat \theta  + z_2 \sin \hat \theta )}  \\
   {g(z_1 \cos \hat \theta  + z_2 \sin \hat \theta )}  \\
\end{array}} \right] = 0.
  \end{equation}

  Premultiplying both sides of (\ref{7}) by $\left[ {\begin{array}{*{20}c}
   { - v_g } & {v_f }  \\
   {u_g } & { - u_f }  \\
\end{array}} \right]^{ - 1} $, we get $\left[ {\begin{array}{*{20}c}
   {f(z_1 \cos \hat \theta  + z_2 \sin \hat \theta )}  \\
   {g(z_1 \cos \hat \theta  + z_2 \sin \hat \theta )}  \\
\end{array}} \right] = 0.$ Then $z_1 \cos \hat \theta  + z_2 \sin \hat \theta =0$,
that is, $z_1$ and $z_2$ are linearly dependent, and
furthermore, $u$ and
  $v$ are linearly dependent.

Here consider the interval $[0,2\pi]$. By \eqref{5}, we have
$T(2\pi)=S(2\pi)=d>0$ and  $T(\pi)=S(\pi)=-d$ where $d$ is shown in \eqref{d}.

Recall the linearly independent vectors
$u=(u_f,u_g),v=(v_f,v_g)\in U$. \eqref{5} together with
that $f$ and $g$ have no common zero point except $0\in R^n$ shows that
\begin{enumerate}
	\item\label{case1'} $S(\theta)=0$ implies $f(\theta)=v_ft$, $g(\theta)=v_gt$
	and $T(\theta)=td$ for some nonzero real number $t$;

	\item\label{case2'} $T(\theta)=0$ implies $f(\theta)=(u_f-v_f)t$, $g(\theta)=(u_g-v_g)t$
	and $S(\theta)=td$ for some nonzero real number $t$;

	\item\label{case3'} $S(\theta)=T(\theta)$ implies $f(\theta)=u_ft$, $g(\theta)=u_gt$
	and $T(\theta)=T(\theta)=td$ for some nonzero real number $t$,
\end{enumerate}
where $t$ cannot be $0$, because there exists no $\theta'$ such that
$f(\theta')=g(\theta')=0$.


Denote the set of zero points of $S$ that are not minima or maxima in the
interval $[0,2\pi]$ by ${\bf0}$.
It is to get $S(\theta)=-S(\theta+\pi)$ for all $\theta\in R$. Then
$|{\bf0}\cap(0,\pi)|=|{\bf0}\cap(\pi,2\pi)|:=l$ is an odd number.
We also have if $S(\theta)=0$, then $T(\theta)T(\theta+\pi)<0$. Hence
$|\{w:w\in{\bf0}\cap(0,2\pi),T(w)>0\}|=|\{w:w\in{\bf0}
\cap(0,2\pi),T(w)<0\}|=l$.
Denote ${\bf0}$ by $\{0_1,\cdots,0_{2l}\}$,
where $0<0_1<\cdots<0_{2l}<2\pi$.

\begin{enumerate}
	\item Assume $l=1$. Based on the foregoing discussion,
		it holds that $T(0_1)T(0_2)<0$.
		Then either $\Lambda([0,0_1]\cap
		\{\theta:T(\theta)>0\})\supset[0,1]$,
		or $\Lambda([0_2,2\pi]\cap\{\theta:T(\theta)>0\})\supset[0,1]$.
	\item Assume $l>1$. If $T(0_{1})>0$, $\Lambda([0,
		0_{1}]\cap\{\theta:T(\theta)>0\})\supset
		[0,1]$. If $T(0_{2l})>0$, $\Lambda([0_{2l},
		\pi]\cap\{\theta:T(\theta)>0\})\supset[0,1]$.
		If  $T(0_{1})<0$ and
		$T(0_{2l})<0$,
		there exists $1\le i<l$ such that $T(0_{2i})
		T(0_{2i+1})<0$. Suppose the contrary: If for each
		$1\le j<l$, $T(0_{2j})T(0_{2j+1})>0$, then
		$|\{w:w\in{\bf0}\cap(0,2\pi),T(w)>0\}|$ is an even number,
		which is a contradiction.
		Since $S(\theta)>0$ for all $\theta\in(0_{2i},0_{2i+1})$,
		$\Lambda([0_{2i},0_{2i+1}]\cap\{\theta:T(\theta)>0\}
		)\supset[0,1]$.
\end{enumerate}

Hence \eqref{judgement} holds.

Based on the above discussion, the set $U$ is a convex set if
$k$ is odd and $f$ and $g$ have no common zero  point except $0\in R^n$.

Given nonzero $(a_1,a_2)\in U$, then $(-a_1,-a_2)
\in U$, and then $\{(a_1t,a_2t):t\in R\}\subset U$. Thus $\{(a_1t,a_2t):t\in
R\}= U$ may hold (see Example \ref{exam1}).  If $\{(a_1t,a_2t):t\in
R\}\ne U$, there exists nonzero $(a_3,a_4)\in U$ such that $a_1a_4-a_2a_3\ne0$,
then $(-a_3,-a_4)\in U$, and then $U=R^2$
by the convexity of the set $U$ and the
homogeneity of $f$ and $g$.
\end{proof}

\subsection{A Counterexample for Switched Polynomial Systems with Three Sub-systems}

In this subsection, we give an example showing that even if a
switched system with three sub-systems admits a LFHD, there may exist no
convex combination of its sub-systems admitting a LFHD.

\begin{example}\label{thm8}
  Let $S$  be a switched linear system with three sub-systems, and the three sub-system matrices are
  \begin{equation*}
  \begin{split}
	  &S_1 :A_1  = \left[ {\begin{array}{*{20}c}
   1 & 0  \\
   0 & { - 1}  \\
\end{array}} \right],\quad S_2 :A_2  = \left[ {\begin{array}{*{20}c}
   { - \sqrt 3 } & { - 1}  \\
   { - 1} & {\sqrt 3 }  \\
\end{array}} \right],\quad S_3 :A_3  = \left[ {\begin{array}{*{20}c}
   { - \sqrt 3 } & 1  \\
   1 & {\sqrt 3 }  \\
\end{array}} \right].
\end{split}\end{equation*}

Firstly, it is easy  to obtain that $V(x) = \frac{1}{2}\left( {x_1^2  + x_2^2 } \right)$ is a LFHD (see Fig. \ref{fig5}).

\begin{figure}
\begin{center}
\includegraphics[height=10cm]{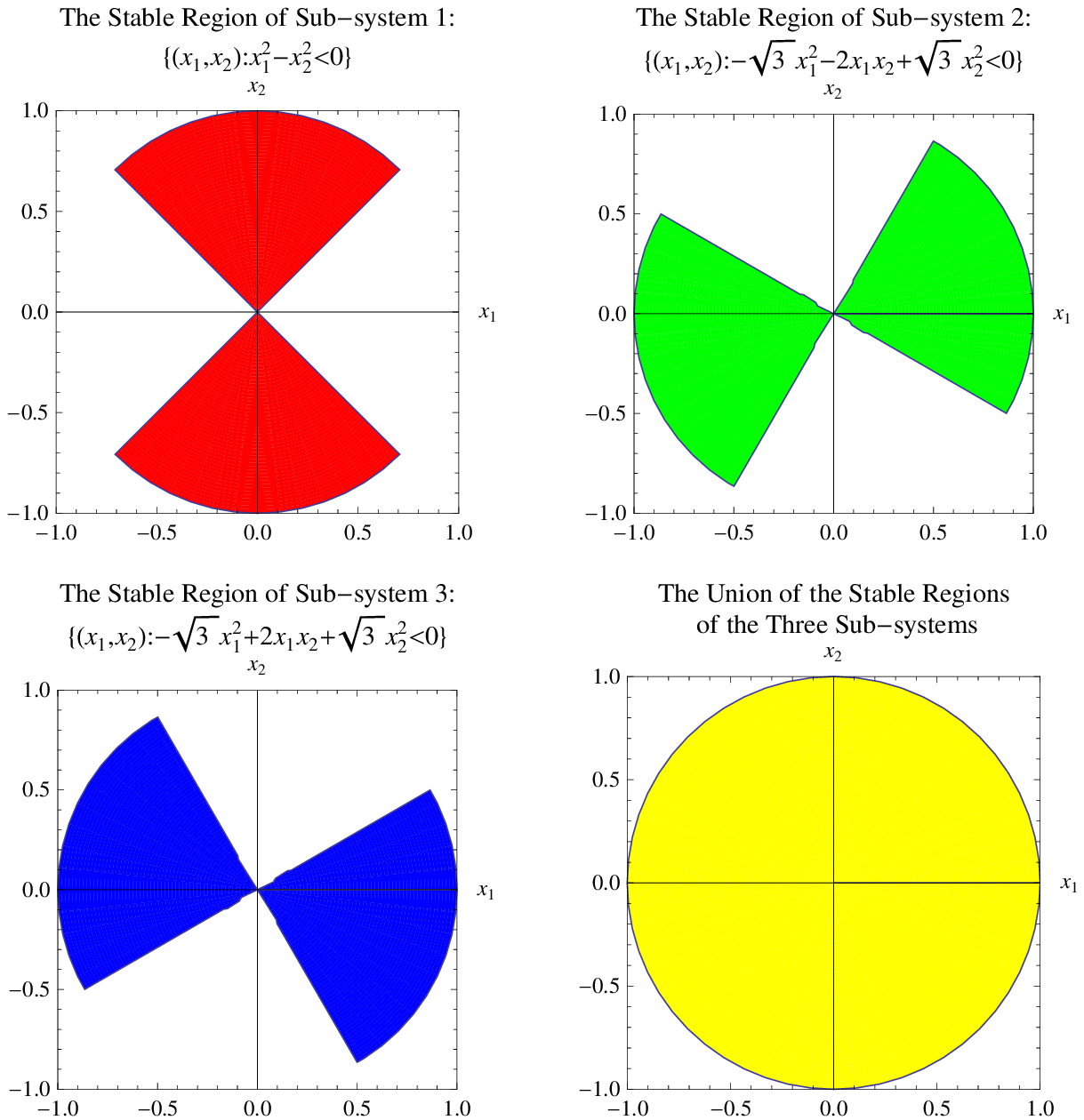}    
\caption{The stable regions of Example \ref{thm8} on the unit disk}  
\label{fig5}                                 
\end{center}                                 
\end{figure}

  Secondly, we prove none of the linear combinations of the three
  sub-systems is asymptotically stable.

  Denote\begin{equation*}\begin{split}
 A &= \sum\limits_{i = 1}^3 {\lambda _i A_i }
  = \left[ {\begin{array}{*{20}c}
   {\lambda _1  - \sqrt 3 \lambda _2  - \sqrt 3 \lambda _3 } & { - \lambda _2  + \lambda _3 }  \\
   { - \lambda _2  + \lambda _3 } & { - \lambda _1  + \sqrt 3 \lambda _2  + \sqrt 3 \lambda _3 }  \\
\end{array}} \right],
\end{split}\end{equation*}
  where $\lambda _1 ,\lambda _2 ,\lambda _3 $ are real variables.
  Notice that $A$  is real symmetric, then we have $A$  either
  is a zero matrix, or has a positive eigenvalue and a negative eigenvalue.
  Hence, system $\dot x = Ax$ is either stable or unstable, but cannot be asymptotically
  stable.

 We can easily get  the unique stable convex combination is $A^*  = \lambda
_1^* A_1  + \lambda _2^* A_2  + \lambda _3^* A_3=0
  $, where $\lambda _1^*  = \frac{{2\sqrt 3 }}{{2\sqrt 3  + 2}}$, $\lambda _2^*  = \lambda _3^*  = \frac{1}{{2\sqrt 3  +
  2}}$.

 Taking $f(x) = - \left.{\dot V(x)} \right|_{S_1 }$,
 $g_i(x)=\left.{\dot V(x)} \right|_{S_{i+1} }$, $i=1,2$, we have the set
$\{(f(x),g_1 (x),\\g_2 (x)):x\in R^n\}$ is not convex.
\end{example}

\subsection{Extended S-Lemma}

In this subsection, based on Theorem \ref{thm3}, by borrowing the idea of
Yakubovich, we give some extended versions of the S-Lemma under some mild
conditions.

In the SHS-Lemma,
if $f$ and $g$ have a nonzero common
zero point, the set $\{(f(x),g(x)):x\in R^n\}$ may be neither
convex nor an angular sector (see Fig. \ref{fig10}).
Luckily, Item $(i)$ of Theorem \ref{thm13} implies $f$ and $g$ have no nonzero
common zero point. So the case shown in Fig. \ref{fig10} does not happen.
In fact, $f$ is strictly copositive with $g$ implies
$f$ is copositive with $g$ and $f$ and $g$ have no nonzero common zero point.
However, $f$ is copositive with $g$ does not
imply $f$ and $g$ have no nonzero common zero point. And $f$ is copositive
with $g$ does not imply there exists $\xi\ge0$ such that $f-\xi g$ is
nonnegative (see Example \ref{exam3}). Later, under the assumption
that the two polynomials
considered have no nonzero common zero point and some extra mild assumptions,
we extend the S-Lemma to the case
of two homogeneous polynomials of the same degree greater than $2$ (Theorem
\ref{thm9}). And based on Theorem \ref{thm9},
we extend the S-Lemma to the case of non-homogeneous polynomials of the same
degree greater than $2$ (Theorem \ref{thm10}) and of the same even degree
greater than $2$ (Theorem \ref{thm11}).

\begin{example}\label{exam3}
	Recall Remark \ref{rem1}. Choose
	$f(x,y)=-x^3+y^3$ and $g(x,y)=y^3-\frac{1}{2}x^3-\frac{1}{2}xy^2$.
	$f$ and $g$ have the common nonzero zero point $(1,1)$. It can be
	calculated that the two boundaries (see Fig. \ref{fig10}) of the set
	$\{(-x^3+y^3,y^3-\frac{1}{2}x^3-\frac{1}{2}xy^2):x,y\in R \}$
	are $y=\frac{1}{2}x$ and  $y=\frac{7}{6}x$. Then $f$ is copositive with
	$g$, but not strictly copositive with $g$. Further we have there exists
	no $\xi\ge0$ such that $f-\xi g$ is nonnegative for all $x,y\in R$.
\end{example}

\subsubsection{Homogeneous  S-Lemma}

\begin{theorem}[HS-Lemma]\label{thm9}
	Let $f,g:R^n\rightarrow R$ be homogeneous polynomials of degree $k\ge1$.
	If $f$ and $g$ have no common zero point except $0\in R^n$ and there
	exists at most one vector $(a,b)\in R^2$ such that $a^2+b^2=1$,
	$a+\delta(a)b>0$ and both
	$\left\{
  \begin{array}{l}
	  f(x)=a\\
	  g(x)=b
  \end{array}
  \right.
  $ and $\left\{
  \begin{array}{l}
	  f(x)=-a\\
	  g(x)=-b
  \end{array}
  \right.
  $ have a solution, where $\delta(t)=\left\{\begin{array}{l}
	  1,\quad\mbox{ if }t=0,\\
	  0,\quad\mbox{ if }t\ne0,
  \end{array}\right.$
	that is to say, there exists at most one straight line passing through
    the origin that is contained in the set $U$,
    then the following two items are equivalent:

	$(i)$ $f$ is copositive with $g$;

	$(ii)$ there exists $\xi\ge0$ such that $f(x)-\xi g(x)\ge0$
	for all $x\in R^n$. In particular, if there exists a vector $z\in R^n$
	such that $f(z)<0$, then there exists $\xi>0$ such that
	$f(x)-\xi g(x)\ge0$ for all $x\in R^n$.
	\label{HS-Lemma}
\end{theorem}

\begin{proof}
	$(ii)\Rightarrow(i)$ holds obviously. We only prove $(i)\Rightarrow(ii)$.

	Next we assume $f$ and $g$ have no common zero point except $0\in R^n$
	and there
	exists at most one vector $(a,b)\in R^2$ such that $a^2+b^2=1$,
	$a+\delta(a)b>0$ and both
	$\left\{
  \begin{array}{l}
	  f(x)=a\\
	  g(x)=b
  \end{array}
  \right.
  $ and $\left\{
  \begin{array}{l}
	  f(x)=-a\\
	  g(x)=-b
  \end{array}
  \right.
  $ have a solution.
	Then by Theorem \ref{thm3}, the set $\{(f(x),g(x)):x\in R^n\}:=U$ is
	a straight line passing through the origin if $k$ is odd, and is
	an angular sector of angle no greater than $\pi$ if $k$ is even.

	Since $f$ is copositive with $g$, $$U
	\cap\{(u,v):u<0,v\ge0\}=\emptyset.$$

	Then there exist $\xi_1<0$ $\xi_2\ge0$ such that $$
	\xi_1f(x)+\xi_2g(x)\le0$$ for all $x\in R^n$.

	Setting $\xi=-\xi_2/\xi_1\ge0$, we have $$
	f(x)-\xi g(x)\ge0$$ for all $x\in R^n$.

	In particular, if there exists $z\in R^n$ such that $f(z)<0$, there
	exist $\xi_1<0$ $\xi_2>0$ such that $$
	\xi_1f(x)+\xi_2g(x)\le0$$ for all $x\in R^n$. Hence there
	exists $\xi>0$ such that
	$$f(x)-\xi g(x)\ge0$$ for all $x\in R^n$.

%
\end{proof}

\subsubsection{Non-homogeneous S-Lemma}

For convenience, we use the semi-tensor product of matrices to represent
a polynomial hereinafter. The concept of the semi-tensor product of matrices is
referred to \cite{Cheng:11} and the references therein. Here we only introduce
the semi-tensor product of two column vectors.

\begin{definition}[\cite{Cheng:11}]\label{def1}
	Let $u,v\in R^n$ be two column vectors, set $u=(u_1,u_2,\cdots,u_n)^T$,
	the semi-tensor product of $u$ and $v$ is defined as $$
	u\ltimes v=(u_1v^T,u_2v^T,\cdots,u_nv^T)^T;$$

	since the semi-tensor product preserves the associative law, $u^m$ is defined
	as $\underbrace{u\ltimes u\ltimes\cdots\ltimes u}_{m}$ inductively.
\end{definition}

By Definition \ref{def1}, a polynomial $f(x):R^n\rightarrow R$ of degree $k$
can be represented as $$f(x)=f_0+f_1x^1+\cdots+f_kx^k,$$
where each $f_i\in R^{n^i}$ is a constant row vector, $i=0,1,\cdots,k$, called
coefficient vector.
Note that $f_0$ and $f_1$ are unique, but other coefficient vectors may be not.

Based on Theorem \ref{thm9}, we have the following Theorems \ref{thm10}
and \ref{thm11}.

\begin{theorem}[NHS-Lemma]\label{thm10}
	Let $f,g:R^n\rightarrow R$ be polynomials of degree $k\ge1$ in the
	form of
	$$f(x)=f_0+f_1x+\cdots +f_kx^k\mbox{ and }
	g(x)=g_0+g_1x+\cdots +g_kx^k,$$
	where $f_i,g_i\in R
	^{n^i}$ are constant row vectors, $i=0,1,\cdots,k$.

	Let us introduce homogeneous functions:
	$$\tilde f:R^{n+1}\rightarrow R,\tilde f(x,t)=f_0t^k+f_1xt^{k-1}+\cdots
	+f_kx^k,$$
	$$\tilde g:R^{n+1}\rightarrow R,\tilde g(x,t)=g_0t^k+g_1xt^{k-1}+\cdots
	+g_kx^k.$$
	Assume
	$\tilde f$ is copositive with $\tilde g$, $\tilde f$ and $\tilde g$
	have no common zero point except $0\in R^{n+1}$, and
	there
	exists at most one vector $(a,b)\in R^2$ such that $a^2+b^2=1$,
	$a+\delta(a)b>0$ and both
	$\left\{
  \begin{array}{l}
	  \tilde f(x,t)=a\\
	  \tilde g(x,t)=b
  \end{array}
  \right.
  $ and $\left\{
  \begin{array}{l}
	  \tilde f(x,t)=-a\\
	  \tilde g(x,t)=-b
  \end{array}
  \right.
  $ have a solution, where $\delta(\cdot)$ is seen in Theorem \ref{thm9}.
	Then there exists $\xi\ge0$ such that $f(x)-\xi g(x)\ge0$
	for all $x\in R^n$.
	\label{S-Lemma}
\end{theorem}

\begin{proof}
	Note that $\tilde f$ is copositive with $\tilde g$ implies $f$ is
	copositive with $g$ by taking $t\equiv1$, and $\tilde f$ and $\tilde g$ have no common
	nonzero zero point implies $f$ and $g$ have no common zero point.
	But the converse is not true.

	Then by Theorem
	\ref{HS-Lemma}, there exists $\xi\ge0$ such that
	$$\tilde f(x,t)-\xi\tilde g(x,t)\ge0\mbox{ for all }(x,t)\in R^{n+1}.$$

	Choosing $t\equiv1$, we have
	$$f(x)-\xi g(x)\ge0\mbox{ for all }x\in R^{n}.$$
\end{proof}

\begin{theorem}\label{thm11}
	Let $f,g:R^n\rightarrow R$ be polynomials of even degree $k$.
	Assume $f$ is copositive with $g$, $f$ and $g$ have no common zero
	point. 
	Denote 	
	\begin{equation*}
		\begin{split}
	&f(x):=f_0+f_1x+\cdots +f_kx^k,\\
	&g(x):=g_0+g_1x+\cdots +g_kx^k,
		\end{split}
	\end{equation*}
	where $f_i,g_i\in R
	^{n^i}$ are constant row vectors, $i=0,1,\cdots,k$.

	Assume that
	homogeneous polynomials $f_ky^k$ and $g_ky^k$ have no common zero point
	except $0\in R^{n}$, $f_ky^k$ is copositive with $g_ky^k$, and assume
	\begin{enumerate}
		\item\label{a1}
	there exist no nonzero vector $(a,b)\in R^2$ such that both
	$\left\{\begin{array}{l}
		f_kx^k=a\\
		g_kx^k=b
	\end{array}\right.$ and\\
	$\left\{\begin{array}{l}
		f_kx^k=-a\\
		g_kx^k=-b
	\end{array}\right.$ have a solution,
		\item\label{a2}
	there exist no $a,b,c,d\in R$ such that $ad-bc=0$, either $ac<0$
	or $bd<0$, and both
	$\left\{\begin{array}{l}
		f(x)=a\\
		g(x)=b
	\end{array}\right.$ and
	$\left\{\begin{array}{l}
		f(x)=c\\
		g(x)=d
	\end{array}\right.$ have a solution,
		\item\label{a3}
	there exist no $a,b,c,d\in R$ such that $ad-bc=0$, either $ac<0$
	or $bd<0$, and both
	$\left\{\begin{array}{l}
		f(x)=a\\
		g(x)=b
	\end{array}\right.$ and
	$\left\{\begin{array}{l}
		f_kx^k=c\\
		g_kx^k=d
	\end{array}\right.$ have a solution,
		\item\label{a4}
	and there exist no $a,b,c,d\in R$ such that $ad-bc=0$, either $ac<0$
	or $bd<0$, and both
	$\left\{\begin{array}{l}
		f_kx^k=a\\
		g_kx^k=b
	\end{array}\right.$ and
	$\left\{\begin{array}{l}
		f(x)=c\\
		g(x)=d
	\end{array}\right.$ have a solution.
	\end{enumerate}

	Then
	there exists $\xi\ge0$ such that $f(x)-\xi g(x)\ge0$ for all
	$x\in R^n$.
\end{theorem}

\begin{proof}
	Let us introduce homogeneous functions:
	$$\tilde f:R^{n+1}\rightarrow R,\tilde f(x,t)=f_0t^k+f_1xt^{k-1}+\cdots
	+f_kx^k$$ and
	$$\tilde g:R^{n+1}\rightarrow R,\tilde g(x,t)=g_0t^k+g_1xt^{k-1}+\cdots
	+g_kx^k.$$

	Firstly we prove $\tilde f$ is copositive with $\tilde g$. That is to say,
	we prove $\tilde f(x,t)<0$ and $\tilde g(x,t)\ge0$ have no common
	solution. Suppose the contrary:
	Assume that there exists $(x_1,t_1)$ such that
	$\tilde f(x_1,t_1)<0$ and $\tilde g(x_1,t_1)\ge0$.
	\begin{enumerate}
	\item
	If $t_1\ne0$, then
	$$f(x_1/t_1)=\tilde f(x_1,t_1)/t_1^k<0,$$
	$$g(x_1/t_1)=\tilde g(x_1,t_1)/t_1^k\ge0,$$ which contradicts that
	$f$ is copositive with $g$.

	\item
	If $t_1=0$, then
	\begin{equation}\label{contra3}
		f_kx_1^k<0\mbox{ and }g_kx_1^k\ge0,
	\end{equation}
	that is to say, $f_ky^k$ is not copositive with $g_ky^k$, which is
	a contradiction.

	\end{enumerate}

	Secondly we prove $\tilde f$ and $\tilde g$ have no common nonzero
	zero point. Suppose the contrary: We assume that  there exists  nonzero
	$(x_2,t_2)\in R^{n+1}$ such that
	$\tilde f(x_2,t_2)=\tilde g(x_2,t_2)=0$.
	\begin{enumerate}
		\item If $t_2\ne0$,
					$$f(x_2/t_2)=\tilde f(x_2,t_2)/t_2^k=0,
					$$ $$g(x_2/t_2)=\tilde g(x_2,t_2)/t_2^k
					=0,$$ which contradicts that $f$ and
					$g$ have no common zero point.
		\item If $t_2=0$, then $x_2\ne0$ and $f_kx_2^k=g_kx_2^k=0$,
			which contradicts $f_ky^k$ and $g_ky^k$ have no common
			zero point except $0\in R^{n}$.

	\end{enumerate}

	Thirdly we prove
	there exist no nonzero vector $(a,b)\in R^2$ such that both
	$\left\{\begin{array}{l}
		\tilde f(x,t)=a\\
		\tilde g(x,t)=b
	\end{array}\right.$ and
	$\left\{\begin{array}{l}
		\tilde f(x,t)=-a\\
		\tilde g(x,t)=-b
	\end{array}\right.$ have a solution. Suppose the contrary:
	If there exist $a,b\in R$,
	$(x_1,t_1),(x_2,t_2)\in R^{n+1}$ such that $a^2+b^2\ne0$,
	$\left\{\begin{array}{l}
		\tilde f(x_1,t_1)=a\\
		\tilde g(x_1,t_1)=b
	\end{array}\right.$ and
	$\left\{\begin{array}{l}
		\tilde f(x_2,t_2)=-a\\
		\tilde g(x_2,t_2)=-b
	\end{array}\right.$, then $(x_1,t_1)\ne0$ and $(x_2,t_2)\ne0$.
	\begin{enumerate}
		\item If $t_1\ne0$ and $t_2\ne0$, then
		\begin{equation*}
			\begin{split}
			&f(x_1/t_1)=\tilde f(x_1,t_1)/t_1^k=a/t_1^k,\quad
			g(x_1/t_1)=\tilde g(x_1,t_1)/t_1^k=b/t_1^k,\\
			&f(x_2/t_2)=\tilde f(x_2,t_2)/t_2^k=-a/t_2^k,\quad
			g(x_2/t_2)=\tilde g(x_2,t_2)/t_2^k=-b/t_2^k,
			\end{split}
		\end{equation*}
		which contradicts Item \ref{a2} in Theorem \ref{thm11}.
		\item If $t_1=0$ and $t_2=0$, then
		\begin{equation*}
			\begin{split}
			f_kx_1^k=a,\quad
			g_kx_1^k=b,\quad
			f_kx_2^k=-a,\quad
			g_kx_2^k=-b,
			\end{split}
		\end{equation*}
		which contradicts Item \ref{a1} in Theorem \ref{thm11}.
		\item If $t_1\ne0$ and $t_2=0$, then
		\begin{equation*}
			\begin{split}
			&f(x_1/t_1)=\tilde f(x_1,t_1)/t_1^k=a/t_1^k,\quad
			g(x_1/t_1)=\tilde g(x_1,t_1)/t_1^k=b/t_1^k,\\
			&f_kx_2^k=-a,\quad
			g_kx_2^k=-b,
			\end{split}
		\end{equation*}
		which contradicts Item \ref{a3} in Theorem \ref{thm11}.
		\item If $t_1=0$ and $t_2\ne0$, then
		\begin{equation*}
			\begin{split}
			&f_kx_1^k=a,\quad
			g_kx_1^k=b,\\
			&f(x_2/t_2)=\tilde f(x_2,t_2)/t_2^k=-a/t_2^k,\quad
			g(x_2/t_2)=\tilde g(x_2,t_2)/t_2^k=-b/t_2^k,
			\end{split}
		\end{equation*}
		which contradicts Item \ref{a4} in Theorem \ref{thm11}.
	\end{enumerate}

	Based on the above discussion, by Theorem \ref{thm9},
	there exists $\xi\ge0$ such that $\tilde f(x,t)-\xi\tilde
	g(x,t)\ge0$ for all $(x,t)\in R^{n+1}$. Taking $t\equiv1$, we have
	$f(x)-\xi g(x)\ge0$ for all $x\in R^{n}$.
\end{proof}

In order to illustrate Theorem \ref{thm11}, we give the following Example
\ref{exam4}.

\begin{example}\label{exam4}
	Consider $f(x_1,x_2)=-7x_1^6 - 2x_1^3 x_2^3 +5x_2^6 - 2x_1^4x_2^2-2$ and
	$g(x_1,x_2)=- 5x_1^6  - x_1^3 x_2^3  + x_2^6 - x_1^4x_2^2-1$ both from
	$R^2$ to $R$.

	$f$ and $g$ are both polynomials of even degree $6$. Easily we have
	$$\left\{
	\begin{array}[]{l}
		-7x_1^6 - 2x_1^3 x_2^3 +5x_2^6 - 2x_1^4x_2^2-2<0\\
		- 5x_1^6  - x_1^3 x_2^3  + x_2^6 - x_1^4x_2^2-1\ge0
	\end{array}
	\right.\Rightarrow -3(x_1^6+x_2^{6})>0.$$ Then
	$f$ is copositive with $g$.
	We easily get that
	$f$ and $g$ have no common zero point, since $$
	\left\{
	\begin{array}[]{l}
		-7x_1^6 - 2x_1^3 x_2^3 +5x_2^6 - 2x_1^4x_2^2-2=0\\
		- 5x_1^6  - x_1^3 x_2^3  + x_2^6 - x_1^4x_2^2-1=0
	\end{array}
	\right.\Rightarrow x_1^6+x_2^{6}=0\Rightarrow x_1=x_2=0.$$
	By \eqref{9} we have
	$-7x_1^6 - 2x_1^3 x_2^3 +5x_2^6 - 2x_1^4x_2^2$ and
	$- 5x_1^6  - x_1^3 x_2^3  + x_2^6 - x_1^4x_2^2$ have no common zero
	point except $(0,0)$, the former is copositive with the latter and
	Item \ref{a1} in Theorem \ref{thm11} holds.
    It is easy to get Items \ref{a2}, \ref{a3} and \ref{a4} in Theorem \ref{thm11}
	hold.

	Based on the above discussion, by Theorem \ref{thm11}, there exists
	$\xi\ge0$ such that $f(x_1,x_2)-\xi g(x_1,x_2)\ge0$ for all $x_1,
	x_2\in R$.

	Since each $\frac{1}{5} < \lambda < \frac{37}{93}$ makes system \eqref{12}
	asymptotically stable, here we might as well choose
	$\lambda=\frac{1}{3}$, i.e., choose $\xi=2$. We have
	$f(x_1,x_2)-2 g(x_1,x_2)=3(x_1^6+x_2^6)\ge0$ for all $x_1, x_2\in R$.

\end{example}

\section{Conclusions}\label{sec3}

This paper studied the relationship between $(i)$ a switched
nonlinear system admits a LFHD
and $(ii)$ the system has a convex combination of its sub-systems
that admits a LFHD. By
using the strict homogeneous S-Lemma presented and proved in this paper,
a necessary and sufficient condition was given under which
$(i)$ is equivalent to $(ii)$ when the system has
two sub-systems, and a counterexample was given to show that
$(i)$ does not imply $(ii)$ when the system has more than two sub-systems.

Besides, the S-Lemma
was extended from quadratic polynomials to polynomials of degree more than
$2$ under some mild conditions.

\section{Acknowledgments}
The authors thank Dr. Ragnar Wallin for supplying some of the references, and
thank the anonymous reviewers and the AE for their valuable comments that
led to an improvement for the readability and  an increase of the range
of applications of the manuscript.

\end{document}